\documentclass[11pt]{amsart}

\usepackage{amsmath,amssymb,amsfonts,amsthm,mathtools}
\usepackage{hyperref}
\usepackage{enumitem}
\usepackage{mathrsfs}

\usepackage{a4wide}

\usepackage[dvipsnames]{xcolor}

\usepackage[normalem]{ulem}

\newtheorem{theorem}{Theorem}[section]
\newtheorem{prop}[theorem]{Proposition}

\newtheorem{lem}[theorem]{Lemma}
\newtheorem{cor}[theorem]{Corollary}
\theoremstyle{remark}
\newtheorem{remark}[theorem]{Remark}
\newtheorem{remarks}[theorem]{Remarks}
\newtheorem{definition}[theorem]{Definition}

\newcommand{\Gal}{\mathrm{Gal}}
\newcommand{\End}{\mathrm{End}}
\newcommand{\Aut}{\mathrm{Aut}}
\newcommand{\GL}{\mathrm{GL}}
\newcommand{\Z}{\mathbf{Z}}
\newcommand{\Q}{\mathbf{Q}}
\newcommand{\F}{\mathbf{F}}

\newcommand{\Qp}{\mathbf{Q}_p}
\newcommand{\Fp}{\mathbf{F}_p}

\newcommand{\Kcyc}{K^{\mathrm{cyc}}}
\newcommand{\Kkum}{K^{\mathrm{Kum}}}

\newcommand{\Kbar}{\overline K}

\title[Cohomology over the maximal Kummer extension]{Cohomology of varieties over the maximal Kummer extension of a number field}

\author{Davide Lombardo and Tam\'as Szamuely}
\date{}

\address{Universit\`a degli studi di Pisa\\
  Pisa, Italy} 
\email{davide.lombardo@unipi.it}  
\email{tamas.szamuely@unipi.it}

\begin{document}

\begin{abstract}
Let $X$ be a smooth projective geometrically connected variety defined over a number field $K$. We prove that the geometric \'etale cohomology of $X$ with $\Q/\Z$-coefficients has finitely many classes invariant under the Galois group of the maximal Kummer extension of $K$ in odd degrees. In particular, every abelian variety has finite torsion over the maximal Kummer extension. This improves results by R\"ossler and the second author as well as Murotani and Ozeki.   We also show that finiteness of torsion of a given abelian variety over non-abelian solvable extensions of $K$ is not controlled by the Galois group of the extension.
\end{abstract}

\maketitle

\section{Introduction}

Let $K$ be a number field with fixed algebraic closure $\overline K$, and let $X$ be a smooth projective geometrically connected variety over $K$. We denote by $\overline{X}$ the base change of $X$ to $\overline{K}$, and for an intermediate field $M$ with $\overline K\supset M\supset K$ we write $G_M:=\Gal(\overline K | M)$. In this note we shall study the question of finiteness of the subgroup of $G_M$-invariants of the étale cohomology groups $H^i_{\textup{ét}}(\overline{X}, \Q/\Z(j))$ over certain infinite extensions $M|K$. Note that when $X=A$ is an abelian variety with dual abelian variety $A^*$, the group $H^1_{\textup{ét}}(\overline{A^*}, \Q/\Z(1))$ identifies with the torsion subgroup of $A(\Kbar)$ as a $G_K$-module (see, for instance, the introduction of \cite{RoesslerSzamuely}), so in this case $G_M$-invariants on cohomology correspond to $M$-rational torsion points of $A$.

In the paper just mentioned, R\"ossler and the second author considered the case when $M$ is obtained by adjoining all roots of unity to $K$ and proved finiteness of $H^i_{\textup{ét}}(\overline{X}, \Q/\Z(j))^{G_M}$ for {\bf odd} $i$. The method of proof generalized that of a much earlier theorem of Ribet \cite{Ribet1981CyclotomicTorsion} who proved finiteness of the torsion subgroup of $M$-points of an abelian variety for $M$ as above.

But there are other infinite extensions $M|K$ over which finiteness of $H^i_{\textup{ét}}(\overline{X}, \Q/\Z(j))^{G_M}$ holds for {odd} $i$. Indeed, recently Murotani and Ozeki \cite{MurotaniOzeki} proved this for a large class of infinite Kummer extensions $M$ of $K$. Here we extend their statement to all Kummer extensions.

\begin{theorem}\label{mainthm}
Let $K^{\rm Kum}$ be the maximal Kummer extension of $K$, obtained by adjoining all roots of elements of $K$. The group  $H^i_{\textup{ét}}(\overline{X}, \Q/\Z(j))^{G_{K^{\operatorname{Kum}}}}$ is finite for all odd $i$ and arbitrary $j$. In particular, the torsion subgroup of $K^{\rm Kum}$-points of an abelian variety over $K$ is finite.   
\end{theorem}

As already pointed out in \cite[Remark 3.6]{RoesslerSzamuely}, for even $i$ the group $H^i_{\textup{ét}}(\overline{X}, \Q/\Z(j))^{G_{K^{\rm Kum}}}$ is infinite because cycle classes on $X_{K^{\rm Kum}}$ generate an infinite subgroup. It is an interesting open question whether the quotient modulo cycle classes is finite.  On the other hand, finiteness may fail already over infinite abelian extensions: for an abelian variety of CM type over $K$ all torsion points are defined over the maximal abelian extension.

The Galois group of the field $K^{\rm Kum}$ over $K$ is  solvable of class 2. The proof of Theorem \ref{mainthm} will proceed by a general reduction to abelian extensions from solvable extensions of finite class: 

\begin{theorem}\label{thm: reduction}
Let $K$ and $X$ be as before, and let $M|K$ be a (possibly infinite) Galois extension contained in $\Kbar$ such that $\Gal(M|K)$ is solvable of finite class. Fix integers $i\geq 0$ and $j$. Assume that for all finite extensions $K' | K$ contained in $\Kbar$ and all subextensions
$
K' \subset F \subset K'M
$
with $F | K'$ abelian the following conditions hold:
\begin{enumerate}
    \item For all primes $p$ we have $\left( H^i_{\textup{ét}}(\overline{X}, \Q_p(j))^{\operatorname{ss}} \right)^{G_F}=0$.
    \item For all but finitely many primes $p$ we have $\left( H^i_{\textup{ét}}(\overline{X}, \Z/p\Z(j))^{\operatorname{ss}} \right)^{G_F}=0$. 
\end{enumerate}
Then $H^i_{\textup{ét}}(\overline{X}, \Q/\Z(j))^{G_M}$ is finite.
\end{theorem}

Here $H^i_{\textup{ét}}(\overline{X}, \Q_p(j))^{\operatorname{ss}}$ (resp.~$H^i_{\textup{ét}}(\overline{X}, \Z/p\Z(j))^{\operatorname{ss}}$) denotes the semisimplification of the module $H^i_{\textup{ét}}(\overline{X}, \Q_p(j))$ (resp.~$H^i_{\textup{ét}}(\overline{X}, \Z/p\Z(j))$) with respect to the action of $\Gal(\Kbar|K)$. It is part of the Tate conjectures (at least with $\Q_p$-coefficients; see \cite{tateconj}, p.~96) that these actions are semisimple and thus the superscripts `ss' can be omitted. If we assume semisimplicity, then in view of well-known arguments about abelian groups (see \cite[Lemma 2.1]{RoesslerSzamuely}) conditions (1) and (2) above can be replaced by the simpler condition that the group $H^i_{\textup{ét}}(\overline{X}, \Q/\Z(j))^{G_F}$ is finite.

\begin{remarks} ${}$
\begin{enumerate}
\item In fact, we prove a slightly stronger statement: there exists a finite extension $K'|K$, depending only on $i$, $j$ and $X$, such that if $M$ is as in the theorem and all abelian extensions $F|K'$ contained in $K'M$ satisfy conditions (1) and (2), then $H^i_{\textup{ét}}(\overline{X}, \Q/\Z(j))^{G_M}$ is finite. On the other hand, there may exist $K''|K$ finite (often we may take $K''=K$; see Remark \ref{rmk: counterexample1} (1)) such that (1) and (2) hold for all abelian $F|K''$ contained in $K''M$ but nevertheless finiteness of $H^i_{\textup{ét}}(\overline{X}, \Q/\Z(j))^{G_M}$ fails. 
\item The semisimplicity conjecture with finite coefficients is perhaps less widely known. However, it is part of general motivic expectations and Corollary 1.3.4 (4) of \cite{Cadoretetal} provides strong evidence for it: Cadoret, Hui and Tamagawa prove that over a finitely generated field of positive characteristic semisimplicity of cohomology with $\Q_p$-coefficients for $p$ large enough (prime to the characteristic) implies semisimplicity with $\Z/p\Z$-coefficients for $p$ large enough. The characteristic 0 analogue of this result would suffice for our purposes. Note also that by an observation of Moonen \cite{Moonen} semisimplicity in characteristic 0 with $\Q_p$-coefficients follows from the usual Tate conjecture on the image of the cycle map.
\end{enumerate}
\end{remarks}

For an abelian variety $A$ over $K$ semisimplicity of the Galois action on the rational Tate module $V_p(A)$ and on $p$-torsion points (for $p$ large enough) is known by classical work of Faltings and Zarhin (see e.g.~Chapter IV by Schappacher in \cite{FaltingsWustholz}). So in this case finiteness of torsion points of $A$ over the maximal abelian extension $(K')^{\rm ab}$ of every finite extension $K' | K$ implies finiteness over solvable extensions of finite class. In particular, plugging in a well-known result of Zarhin (\cite{ZarhinCM}, Corollary of Theorem 1) we obtain: 

\begin{cor}\label{cor: finiteness over solvable extensions for non-CM AV}
Let $A$ be an abelian variety over $K$ having no simple isogeny factor of CM type over the algebraic closure $\Kbar$. Then $A$ has finitely many torsion points over Galois extensions $M|K$ such that $\Gal(M| K)$ is solvable of finite class. \end{cor} This statement was also obtained using arguments similar to ours in very recent work of Huryn \cite{Huryn}. Huryn's preprint appeared while we were writing up the present note and we took the opportunity to borrow one of his arguments to simplify the proof of one of our lemmas (Lemma \ref{lemma: p-adic nilpotent implies algebraic nilpotent} below).

Observe that the conditions on the extensions of $K$ appearing in Theorem \ref{thm: reduction} are purely Galois-theoretic. This prompts the question whether for fixed $X$, $i$ and $j$ there exists a Galois-theoretic criterion on extensions $M|K$ that ensures the finiteness of  $H^i_{\textup{ét}}(\overline{X}, \Q/\Z(j))^{G_M}$. Unfortunately, this is not the case:

\begin{prop}\label{prop:solvable_example} There exist two infinite Galois extensions $M, M'$ of $\Q$ with $\Gal(M | \Q)\cong\Gal(M' | \Q)$ solvable of class 2 such that all abelian varieties defined over subfields of $M$ have finite torsion subgroup over $M$ but some abelian varieties defined over $\Q$ have infinite torsion over $M'$. 
\end{prop}

{After some preliminaries, Theorem \ref{thm: reduction} will be proven in Section \ref{secred}. Theorem \ref{mainthm} will be deduced from Theorem \ref{thm: reduction} in Section \ref{secmainthm} by means of an adjustment of arguments from \cite{RoesslerSzamuely}. The last section contains the proof of Proposition \ref{prop:solvable_example}.
  }
  
  We are very grateful to the referee and the handling editor for their pertinent comments; in particular, the referee has saved us from an embarrassing misuse of Caruso's theorem in the first version.

\section{On a theorem of Serre and Wintenberger}

During the proof of Theorem \ref{thm: reduction} we shall need a slight variation on a theorem contained in Wintenberger's paper \cite{Wintenberger} that has its origin in results explained by Serre in courses given at Collège de France in 1985/86. We begin with a list of conditions coming from Section 3.3 of \cite{Wintenberger}. 

\begin{definition}\label{def: restricted representations}
Let $K$ be a number field, $p$ a prime, $M_p$ a finite-dimensional $\F_p$-vector space, and $\rho_{p} : \operatorname{Gal}(\overline{K} | K) \to \operatorname{GL}({M_p})$ a {continuous} representation. Fix natural numbers $d_0, b$. We say that $\rho_p$ is an \emph{SW representation of type $(d_0, b)$} if the following hold:
    \begin{enumerate}
        \item The $\F_p$-dimension of $M_p$ is at most $d_0$.
        \item The representation $\rho_p$ is semisimple.
        \item Every inertia subgroup at a place of $K$ not lying above $p$ acts on $M_p$ via a pro-$p$-quotient.
        \item The weights of the action of the tame inertia subgroup at every place lying above $p$ on $M_p$ lie in the interval $[0, b]$.
    \end{enumerate}
\end{definition}

In condition (4), tame inertia weights are meant in the sense of Serre's paper \cite{MR387283}. Namely, it is known (see \cite[\S 1.6]{MR387283}) that after restricting to a Jordan--H\"older factor $W_p$ of $M_p$, for a place $v$ lying above $p$ the inertia action factors via the tame inertia quotient $I_v^t$, and is given by a character $I_v^t\to \F_{p^{d}}^\times$, where  $d$ is the $\F_p$-dimension of $W_p$. This character in turn is a product of powers of Serre's fundamental characters $\theta_0,\dots,\theta_{d-1}$ of level $d$; the latter are defined by composing the natural projection $I_v^t\to \mu_{p^d-1}\stackrel\sim\to\F_{p^d}^{\times}$ with each of the $d$ possible field automorphisms of $\F_{p^d}$. The exponents of the $\theta_i$ in the product description are the tame inertia weights; they are well defined up to cyclic reordering \cite[p.~7]{Wintenberger}. This is the convention used in \cite[\S 3.3]{Wintenberger}.

\begin{theorem}\label{thm: Serre-Wintenberger}
    Fix a number field $K$ and natural numbers $d_0$ and $b$. There exists a finite extension $K'|K$  such that for all but finitely many primes $p$ the derived subgroup of $\rho_p\left( \operatorname{Gal}(\overline{K} | K') \right)$ is a perfect group for every SW representation  $\rho_p : \operatorname{Gal}(\overline{K} | K) \to \operatorname{GL}({M_p})$ of type $(d_0, b)$.
\end{theorem}

For the proof we need some auxiliary statements about algebraic groups over finite fields.

\begin{lem}[Borel-Tits]\label{lemma:Borel-Tits}
Let $p$ be a prime number and $N$ be a {connected} semisimple algebraic group over $\F_p$. Let $f : \widetilde{N} \to N$ be the universal cover of $N$ and denote by $N(\F_p)_u$ the image of $\widetilde{N}(\F_p) \xrightarrow{f} N(\F_p)$.
If $p>3$, the group $N(\F_p)_u$ is perfect and coincides with the derived subgroup of $N(\F_p)$.    
\end{lem}
\begin{proof}
    See \cite[6.5 and 6.6]{BorelTits} and \cite[\S 1.2]{Wintenberger}. 
\end{proof}

\begin{prop}\label{prop: derived subgroup}
Let $G$ be a connected reductive group over a finite prime field $\F_p$, and let $N=G'$ be its connected derived subgroup. If $p>3$, we have an equality of groups of $\F_p$-points
\[
G(\F_p)' = N(\F_p)'.
\]
\end{prop}

\newcommand{\ad}{\operatorname{ad}}

\begin{proof} It suffices to prove the containment $G(\F_p)' \subseteq N(\F_p)'$, the other one being obvious. By the lemma above we may replace $N(\F_p)'$ by  $N(\F_p)_u = f\big(\widetilde{N}(\F_p)\big)$.

Write $Z_G$ for the (scheme-theoretic) center of $G$ and $G^{\ad}:=G/Z_G$ for the adjoint quotient.
The commutator morphism
\[
c:G\times G\longrightarrow G,\qquad (g,h)\longmapsto [g,h]:=g^{-1}h^{-1}gh
\]
takes values in $N$ (by definition of the derived subgroup) and is invariant under translation by $Z_G\times Z_G$ on the source. Hence $c$ factors uniquely through the quotient
\[
G\times G\longrightarrow (G/Z_G)\times (G/Z_G)=G^{\ad}\times G^{\ad},
\]
yielding a morphism of $\F_p$-schemes $\alpha:G^{\ad}\times G^{\ad}\to N$ with $c=\alpha\circ\pi$, where $\pi:G\times G\to G^{\ad}\times G^{\ad}$ is the natural map.

Now note that $G$ and $\tilde{N}$ have the same adjoint quotient, since both are given by $G/Z_G \cong N/(Z_G \cap N) \cong \widetilde{N}/Z_{\widetilde{N}}$, where $Z_{\widetilde{N}}$ denotes the center of $\widetilde{N}$. Moreover, since $\widetilde{N}$ is semisimple, it agrees with its own derived subgroup. Thus, the same construction as above, applied to $\widetilde{N}$, produces a morphism
\[
\beta:\widetilde{N}^{\ad}\times \widetilde{N}^{\ad} = G^{\ad}\times G^{\ad}\longrightarrow \widetilde{N}
\]
such that $f\circ\beta=\alpha$. Therefore, at the level of $\F_p$-points we obtain a map
\[
\beta_p:G(\F_p)\times G(\F_p)\longrightarrow \widetilde{N}(\F_p)
\]
such that for $g,h\in G(\F_p)$ we have
\[
[g,h]=f\left(\beta_p(g,h)\right)\in f\big(\widetilde{N}(\F_p)\big),
\]
which concludes the proof.\end{proof}

\begin{proof}[Proof of Theorem \ref{thm: Serre-Wintenberger}]
Given an SW representation $\rho_p$ as above, Serre and Wintenberger construct algebraic subgroups $G_p^{\operatorname{alg}}$ and $N_p$ of $\GL({M_p})$ and prove the existence of a finite extension $K' | K$ such that for all but finitely many $p$  the following hold (see \cite[\S 3.3, Théorème 4]{Wintenberger}):
    \begin{enumerate}
        \item $G_p^{\operatorname{alg}}$ is a reductive group and $N_p$ is its derived subgroup.
        \item We have the containments $$N_p(\F_p)'\subset\rho_p(\operatorname{Gal}(\overline{K} | K'))\subset G_p^{\operatorname{alg}}(\F_p).$$  
    \end{enumerate}
(Note that Wintenberger states this with the subgroup $N_p(\F_p)_u$ instead of $N_p(\F_p)'$ but, as remarked above, the two are the same for $p>3$.)

We now show that $N_p(\F_p)'$ equals the derived subgroup of $\rho_p(\operatorname{Gal}(\overline{K} | K'))$. This will prove the theorem as the group $N_p(\F_p)'$ is perfect by Lemma \ref{lemma:Borel-Tits}.   Using the second containment in (2) we obtain
    \[
    \left[ \rho_p(\operatorname{Gal}(\overline{K} | K')), \; \rho_p(\operatorname{Gal}(\overline{K} | K')) \right] \subseteq \left[ G_p^{\operatorname{alg}}(\F_p), G_p^{\operatorname{alg}}(\F_p)\right] = N_p(\F_p)'
    \]
    where the last equality holds by Proposition \ref{prop: derived subgroup}.
    On the other hand, using the first containment in (2) we have 
    \[
    N_p(\F_p)' = \left[ N_p(\F_p)', \; N_p(\F_p)' \right] \subseteq \left[ \rho_p(\operatorname{Gal}(\overline{K} | K')), \; \rho_p(\operatorname{Gal}(\overline{K} | K'))\right],
    \]
    where the first equality uses that $N_p(\F_p)'$ is perfect. Combining the above inclusions we deduce that the derived subgroup of $\rho_p(\operatorname{Gal}(\overline{K} | K'))$ is the perfect group $N_p(\F_p)'$, as desired.
\end{proof}

\section{Proof of Theorem \ref{thm: reduction}}\label{secred}

We begin the proof of Theorem \ref{thm: reduction}. As already noted in the introduction, the finiteness of $H^i_{\textup{ét}}(\overline{X}, \Q/\Z(j))^{G_M}$ is equivalent to  the vanishing of $H^i_{\textup{ét}}(\overline{X}, \Q_p(j))^{G_M}$ for all primes $p$ and of $H^i_{\textup{ét}}(\overline{X}, \Z/p\Z(j))^{G_M}$ for all but finitely many $p$ (see \cite{RoesslerSzamuely}, Lemma 2.1). Hence we first prove two auxiliary results, one with $\Q_p$-coefficients and one with $\Z/p\Z$-coefficients.

We formulate these for arbitrary finite-dimensional $p$-adic $G_K$-representations $V_p$ stabilizing a $\Z_p$-lattice $T_p\subset V_p$. As before, $V_p^{\rm ss}$ denotes the semisimplification of $V_p$. 

\begin{prop}\label{prop:redQ_p}
 Assume that the Zariski closure of the image of $G_K$ in
$\Aut_{\Qp}(V_p^{\rm ss})$ is connected. Let $M|K$ be a subextension of $\Kbar|K$ with solvable Galois group of finite class. If for all abelian subextensions $F|K$ of $M|K$ we have $\left(V_p^{\rm ss}\right)^{G_F}=0$, 
then $V_p^{G_M}=0$.
\end{prop}

The proof uses a group-theoretic lemma.

\begin{lem}\label{lemma: p-adic nilpotent implies algebraic nilpotent}
Let $G \subset\GL_{n}(\Q_p)$ be a connected reductive algebraic subgroup and $H \subset G(\Q_p)$ be an abstract subgroup, dense in the Zariski topology. Suppose that $H$ is solvable of finite class as an abstract group. Then $G$ is a torus.
\end{lem}

\begin{proof}
    The following proof is inspired by \cite[Proof of Theorem 1.1(a')]{Huryn}. Let $n \geq 1$ be such that $H$ is $n$-step solvable.
    Consider the $n$-fold commutator map
    \[
    D_n :G^{\times 2^n} \to G',
    \]
    where $G'$ denotes the derived subgroup of $G$. We claim that this map is surjective (in the algebraic sense). Indeed, the semisimple group $G'$ coincides with its derived subgroup, which in particular gives surjectivity of $D_1 : G' \times G' \to G'$. Iterating this observation implies that $D_n : (G')^{\times 2^n} \to G'$ is surjective, and hence so is $D_n : G^{\times 2^n} \to G'$.
    By assumption, $H^{\times 2^n}$ is Zariski-dense in $G^{\times 2^n}$, and $D_n$ is trivial on $H^{\times 2^n}$ since $H$ is $n$-step solvable. This implies that $D_n$ is the trivial map. Since it is also surjective, $G'$ is the trivial group. A reductive group with trivial derived subgroup is a torus.
\end{proof}

\begin{proof}[Proof of Proposition \ref{prop:redQ_p}] Assume  $V_p^{G_M}\neq 0$, and hence also $W := (V_p^{\operatorname{ss}})^{G_M}\neq 0$.
As $G_M$ is normal in $G_K$, the $\Qp$-submodule $W\subset V_p^{\operatorname{ss}}$ is $G_K$-stable, and the induced representation
$
\rho_W\colon G_K \longrightarrow \Aut_{\Qp}(W)
$
factors through $\Gal(M| K)$. In particular, $\rho_W(G_K)$ is a solvable group of finite class.

Denote by $G_{W}$ the Zariski closure of $\rho_W(G_K)$ seen as an algebraic group over $\Qp$. It is connected by assumption, and by construction it is a reductive group whose $\Q_p$-points contain $\rho_W(G_K)$ as a Zariski dense subgroup. Lemma \ref{lemma: p-adic nilpotent implies algebraic nilpotent} then shows that $G_W$ is commutative.
The fixed field $F$ of $\ker(G_K \to {G_W(\Qp)})$ is an abelian extension of $K$ with $W\subset\left(V_p^{\operatorname{ss}}\right)^{G_F}$ by construction. This proves $\left(V_p^{\rm ss}\right)^{G_F}\neq 0$. 
\end{proof}

\begin{prop}\label{prop:redmodp}
Assume that the image of $G_K$ in the automorphism group of $(T_p/pT_p)^{\operatorname{ss}}$ has perfect derived subgroup. Let $M|K$ be a pro-solvable subextension of $\Kbar|K$. If for all abelian subextensions $F|K$ of $M|K$ we have $\left((T_p/pT_p)^{\rm ss}\right)^{G_F}=0$, 
then $(T_p/pT_p)^{G_M}=0$.
\end{prop}

\begin{proof} Assume $(T_p/pT_p)^{G_M}\neq 0$, and hence also $W_p := \left( \left(T_p/pT_p \right)^{\operatorname{ss}} \right)^{G_M}\neq 0$.
As in the previous proof, the subspace $W_p$ is $G_K$-stable and the induced representation
$
\rho_{W_p}\colon G_K \longrightarrow \Aut_{\Fp}(W_p)
$
factors through $G=\Gal(M| K)$. Thus $\rho_{W_p}(G_K)$ is a finite solvable group.

Let $\rho_{p}\colon G_K\to\Aut_{\Fp}( (T_p/pT_p)^{\operatorname{ss}} )$ be the full mod $p$ representation on $(T_p/pT_p)^{\operatorname{ss}}$, and let  $
e_{W_p}:\rho_p(G_K)\longrightarrow \Aut_{\F_p}(W_p)$ be the homomorphism induced by restriction to the $G_K$-stable subspace
$W_p$.
Then
\[
\rho_{W_p}(G_K)=e_{W_p}\big(\rho_{p}(G_K)\big)\ \supseteq\ e_{W_p}\big(\rho_{p}(G_K)'\big).
\]
By assumption
$\rho_{p}(G_K)'$ is a perfect group, hence so is $e_{W_p}\big(\rho_{p}(G_K)'\big)$. As the latter group is contained in the solvable group $\rho_{W_p}(G_K)$, it must be trivial. Since applying $e_{W_p}$ commutes with taking the derived subgroup, we get triviality of $e_{W_p}\big(\rho_{p}(G_K)\big)'$, which means that 
$\rho_{W_p}(G_K)=e_{W_p}\big(\rho_{p}(G_K)\big)$ is abelian. The fixed field $F$ of $\ker(\rho_{W_p})$ is then an abelian extension of $K$ contained in $M$ such that
$W_p\subset\left((T_p/pT_p)^{\operatorname{ss}}\right)^{G_F}$. It follows that $\left((T_p/pT_p)^{\operatorname{ss}}\right)^{G_F}\neq 0.$
\end{proof}

\begin{proof}[Proof of Theorem \ref{thm: reduction}] We take up the notation of the theorem. Recall from the beginning of this section that we have to prove $H^i_{\textup{ét}}(\overline{X}, \Q_p(j))^{G_M}=0$ for all primes $p$ and $H^i_{\textup{ét}}(\overline{X}, \Z/p\Z(j))^{G_M}=0$ for all but finitely many $p$. 
Set $V_p:=H^i_{\textup{ét}}(\overline{X}, \Q_p(j))$ and let $T_p$ be $H^i_{\textup{ét}}(\overline{X}, \Z_p(j))$ modulo its torsion subgroup. By a theorem of Larsen and Pink (see \cite{MR1441234} or \cite[Proposition 6.14]{MR1150604}) there exists a finite extension $K'|K$ such that the Zariski closure of the image of $G_{K'}$ in $\Aut_{\Qp}(V_p^{\rm ss})$ is connected. Hence condition (1) of the theorem implies the assumptions of Proposition \ref{prop:redQ_p} with $K'$ in place of $K$, so that $V_p^{G_{K'M}}=0$ and hence also $V_p^{G_M}=0$.

The other half of the proof amounts to verifying that there exists a finite extension $K'|K$ independent of $p$ so that the image of $G_{K'}$ in the automorphism group of $M_p:=(T_p/pT_p)^{\operatorname{ss}}$ has perfect derived subgroup for $p$ large enough, for then applying Proposition \ref{prop:redmodp} with $K'$ in place of $K$ will yield $(T_p/pT_p)^{G_{K'M}}=(T_p/pT_p)^{G_M}=0$ for $p$ large enough. To justify the above claim we are free to twist $M_p$ by a power of the mod $p$ cyclotomic character $\bar\chi_p$; indeed, twisting by $\bar\chi_p$ multiplies all operators by scalars, hence does not affect commutators. So we may and do assume that in the cohomology group $H^i_{\textup{ét}}(\overline{X}, \Z_p(j))$ defining $T_p$ we have $j\geq i$. Under this assumption we shall prove that for $p$ large enough $M_p$ is an SW representation of $G_{K'}$ of type $(d_0,b)$ as in Definition \ref{def: restricted representations} for a suitable {pair $(d_0, b)$ and a suitable} finite extension $K'|K$ independent of $p$. From this the claim will follow by an application of Theorem \ref{thm: Serre-Wintenberger} after replacing $K'$ by another finite extension if needed.

Note that by comparison with singular cohomology over $\bf C$  one knows that the groups $H^i_{\textup{ét}}(\overline{X}, \Z_p(j))$ and $H^{i+1}_{\textup{ét}}(\overline{X}, \Z_p(j))$   are torsion free for $p$ large enough (see e.g.~\cite{RoesslerSzamuely}, Fact 2.3 for references), from which one deduces an isomorphism $T_p/pT_p\cong H^i_{\textup{ét}}(\overline{X}, \Z/p\Z(j))$ using the exact sequence
$$
0\to H^i_{\textup{ét}}(\overline{X}, \Z_p(j))/pH^i_{\textup{ét}}(\overline{X}, \Z_p(j))\to H^i_{\textup{ét}}(\overline{X}, \Z/p\Z(j))\to H^{i+1}_{\textup{ét}}(\overline{X}, \Z_p(j))[p]\to 0. $$

 The same comparison shows that condition (1) of SW representations is satisfied for $M_p$ for $p$ large enough, with $d_0$ the $i$-th Betti number of $X({\bf C})$; { in fact, the $\F_p$-dimension of $M_p$ equals this $d_0$ for $p$ large enough.} Condition (2) being automatic, we pass to (3). As explained in (\cite{RoesslerSzamuely}, Lemma 3.3 (b)), a strong form of Grothendieck's local monodromy theorem implies that there is a finite extension $K'|K$ independent of $p$ so that all inertia subgroups in $G_{K'}$ at places not lying above $p$ act unipotently on $V_p$. It follows that these inertia subgroups act on $T_p/pT_p$ with eigenvalues congruent to 1 mod $p$, and hence they act on $M_p$ via a $p$-group. Finally, condition (4) with $b=j$ is a consequence of Serre's tame inertia conjecture as proven by Caruso \cite{Caruso}. It is here that we use that $j \geq i$: \cite[Theorem 1.2]{Caruso} shows that, for $p$ unramified in $K$ (so for all but finitely many $p$), the tame inertia weights on the $\F_p$-linear dual $H^i_{\textup{ét}}(\overline{X}, \Z/p\Z)^\vee$ lie in the interval $[0, i]$, which for $j \geq i$ implies that the weights of $H^i_{\textup{ét}}(\overline{X}, \Z/p\Z(j))$ lie in $[j-i, j] \subseteq [0, j]$.
\end{proof}

\begin{remarks}\label{rmk: counterexample1} ${}$
\begin{enumerate}
\item One may ask whether Theorem \ref{thm: reduction} holds with $K'=K$, i.e.~whether it is sufficient to check finiteness over abelian subextensions of $M|K$. The answer is negative, even when $i=j=1$ and $X$ is an abelian variety,
as the following example shows. 

Recall from 
the introduction that this case amounts to understanding the torsion subgroups of abelian varieties.
Consider the elliptic curve $E$ with Weierstrass equation $y^2=x^3+x$, with complex multiplication by $\mathbb{Z}[i]$ over $\Q(i)$. By the theory of complex multiplication the field $M$ obtained by adjoining all torsion points of $E$ to $\Q(i)$ is the maximal abelian extension of $\Q(i)$ (it is easy to show that $M| \Q(i)$ is abelian, and this is all we need). This implies that $M| \Q$ is a solvable extension, because both $M| \Q(i)$ and $\Q(i)| \Q$ are abelian. By construction, the torsion subgroup $E(M)_{\mathrm{tors}}\subset E(M)$ is infinite. On the other hand, since $\Q^{\operatorname{ab}}=\Q^{\operatorname{cyc}}$ by the Kronecker--Weber theorem,  Ribet's result \cite{Ribet1981CyclotomicTorsion} shows that $E(\Q^{\operatorname{ab}})_{\mathrm{tors}}$ is finite, hence there is no abelian subextension $F$ of $M| \Q$ such that $E(F)_{\mathrm{tors}}$ is infinite.
\item Note that Proposition \ref{prop:redmodp} holds for \textit{pro}-solvable Galois groups and not just for solvable groups of finite class. Therefore the part of Theorem \ref{thm: reduction} concerning the triviality of $H^i_{\textup{ét}}(\overline{X}, \Z/p\Z(j))^{G_M}$ for all but finitely many $p$ holds for $\Gal(M|K)$ prosolvable. In the case of the $p$-adic Tate module of an abelian variety, combining this more general form with Zarhin's theorem as in the proof of Corollary~\ref{cor: finiteness over solvable extensions for non-CM AV}, one recovers \cite[Theorem~1.1(b)]{Huryn}. We thank Jake Huryn for drawing our attention to this point.
\end{enumerate}
\end{remarks}

\section{Some lemmas from Galois theory}

In this section we collect some statements about Galois extensions of number fields that will be used in the proof of Theorem \ref{mainthm}.

\begin{prop}\label{prop: abelianisation Kummer}
Let $K$ be a number field and denote by $\Kkum$, $K^{\rm ab}$, and $\Kcyc$ its maximal Kummer, abelian and cyclotomic extensions, respectively (inside a fixed algebraic closure). Denote moreover by $w$ the number of roots of unity in $K$ and consider $L:=K\left((K^\times)^{1/w}\right)$, the Kummer extension of $K$ obtained by adjoining all $w$-th roots of elements of $K$.
We then have
$$
\Kkum \cap K^{\mathrm{ab}}=L\Kcyc.
$$
In particular, $\Gal((\Kkum \cap K^{\mathrm{ab}})|\Kcyc)$ is an abelian group of exponent dividing $w$. 
\end{prop}

The proof will use (part of) the following lemma due to Schinzel.

\begin{lem}\label{schinzel}
Let $K$ be a field, $m$ a positive integer prime to the characteristic of $K$, and $a \in K^\times$. The splitting field of the polynomial $x^m-a$ over $K$ is an abelian extension of $K$ if and only if $a^u\in K^{\times m}$, where $u$ is the number of $m$-th roots of unity contained in $K$.    
\end{lem}

\begin{proof}
The original reference is \cite[Theorem 2]{MR429819}; a simpler proof is given in \cite[Theorem 4.1]{stevenhagen}.    
\end{proof} 

\begin{proof}[Proof of Proposition \ref{prop: abelianisation Kummer}] Since $\mu_w\subset K^\times$ by assumption, $L|K$ is an abelian Kummer extension and as such contained in $\Kkum\cap K^{\rm ab}$. For the reverse inclusion, we show that every finite extension $E|\Kcyc$ contained in $\Kkum\cap K^{\rm ab}$ is contained in $L\Kcyc$. Since $E\subset\Kkum$, there exists an integer $m\geq 1$ such that $E\subset \Kcyc\bigl((K^\times)^{1/m}\bigr)$. Since moreover $\mu_m\subset \Kcyc$, Kummer theory shows that there exist $\alpha_1,\dots, \alpha_r\in \Kkum\cap K^{\rm ab}$ such that ${a_i := \alpha_i^m}\in K^\times$ for $i=1,\dots, r$ and $E\subset \Kcyc(\alpha_1,\dots, \alpha_r)$. It thus suffices to show $\alpha_i\in L\Kcyc$ for all $i$. By construction, the splitting field $K(\mu_m,\alpha_i)$ of the polynomial $x^m-a_i$ over $K$ is contained in  $K^{\rm ab}$, so Lemma~\ref{schinzel} applies and yields  $b_i\in K^\times$ such that $a_i^u=b_i^m$, with $u$ the number of $m$-th roots of unity in $K$. It follows that $\left({\alpha_i^u}/{b_i}\right)^m={a_i^u}/{b_i^m}=1$, so $\alpha_i^u/b_i$ is a root of unity of order dividing $m$. Since $\Kcyc$ contains all roots of unity, there exists a root of unity $\zeta_i\in \Kcyc$ such that $\zeta_i^u={\alpha_i^u}/{b_i}$, which can be rewritten as $\alpha_i^u=\zeta_i^ub_i$. This shows $\alpha_i\in\Kcyc \bigl((K^\times)^{1/u}\bigr)$. But $u$ divides $w$, the number of all roots of unity in $K$, and therefore $\alpha_i\in\Kcyc \bigl((K^\times)^{1/w}\bigr) =L
\Kcyc$, as claimed. 
\end{proof}

\begin{definition}
 We shall say that an algebraic extension $L|\mathbf{Q}$ has \emph{bounded local degrees} if there exists
a constant $B\ge 1$ such that for every rational prime $p$ and every place $w\mid p$ of $L$
the completion $L_w$ satisfies $[L_w:\mathbf{Q}_p]\le B$.   
\end{definition}

{Note that a finite extension $K|\Q$ has bounded local degrees (tautologically), and so does any abelian extension $L|K$ of finite exponent, by local class field theory and the structure of the multiplicative group of a $p$-adic field.}

\begin{lem}\label{lem:bounded local degrees}
Let $L|\mathbf{Q}$ be an algebraic extension with bounded local degrees. Let $L^{\rm cyc}$ be its maximal cyclotomic extension, and let $L'\subset L^{\rm cyc}$ be the maximal subextension of $L^{\rm cyc}|L$ which is unramified at all finite places of $L$. Then $L'$ has bounded local degrees.
\end{lem}

\begin{proof}Denote by $B$ a bound for the local degrees of $L|\Q$ as in the above definition.
Consider first a prime number $p$ and a place $w$ of $L$ dividing $p$. 
Let $J\subset \Gal(\Q_p(\mu_{p^\infty})|\Q_p)\simeq \Z_p^\times$ be the image of the inertia subgroup of $\Gal(L_w(\mu_{p^\infty})|L_w)$ under the natural (injective) restriction map $\Gal(L_w(\mu_{p^\infty})|L_w) \to \Gal(\Q_p(\mu_{p^\infty})|\Q_p)$.  We show that $\Z_p^\times/J$ is finite of order at most $B$, for which it suffices to establish the same property for every finite quotient of $\Z_p^\times/J$.  Such a finite quotient corresponds to a finite subextension $F|\Q_p$ of $\Q_p(\mu_{p^\infty})|\Q_p$  contained in the fixed field of $J$. Here $F|\Q_p$ is totally ramified (because so is $\Q_p(\mu_{p^\infty})|\Q_p)$), and
$L_wF|L_w$ is unramified, by definition of $J$. Hence $[F:\Q_p]=e(F|\Q_p)$ divides $e(L_wF|\Q_p)=e(L_wF|L_w)e(L_w|\Q_p)=e(L_w|\Q_p)$. Thus $[F:\Q_p]\leq e(L_w|\Q_p)\leq [L_w:\Q_p]\leq B$, as claimed.
Writing $\Gamma:=\operatorname{Gal}(L^{\rm cyc}| L)$ and viewing $\Gamma$  as a closed subgroup of $\operatorname{Gal}(\Q^{\rm cyc}| \Q)$, we obtain that the inertia subgroup $I_w\subset \Gamma$ maps onto an open subgroup of index $\leq B$ via the $p$-adic cyclotomic character $\Gamma\ \to\ \mathbf{Z}_p^\times$. 

We now identify $\operatorname{Gal}(\Q^{\rm cyc}| \Q)$ with $\hat{\mathbf{Z}}^\times\cong \prod_p\Z_p^\times$ via the adelic cyclotomic character. The Galois group $\operatorname{Gal}(L^{\rm cyc} | L')$ is the closed subgroup of $\Gamma$
generated by the inertia subgroups at all finite places of $L$. Among these, the inertia subgroups of primes dividing $p$ have trivial image in the  $\Z_\ell^\times$-components for $\ell \neq p$, so their images in $\hat{\mathbf{Z}}^\times$ generate a closed subgroup $I_p\subset \Z_p^\times$. By the previous paragraph, we have $[\Z_p^\times:I_p]\leq B$. 
Now set $I:=\prod_p I_p\ \le\ \prod_p \mathbf{Z}_p^\times\simeq \widehat{\mathbf{Z}}^\times$. By the above, the group $\Gal(L' | L)$ identifies with a  subquotient of $\widehat{\mathbf{Z}}^\times/I\cong \prod_p ({\mathbf{Z}}^\times_p /I_p)$ and as such has exponent dividing $B!$ by the arguments above. Since moreover $L'|L$ is an unramified extension, all decomposition subgroups in $\Gal(L'|L)$ above places of $L$ are procyclic and hence must be finite of order dividing $B!$. This shows the boundedness of the local degrees of $L'$ relative to $L$, from which the lemma follows.   
\end{proof}

We note for later use the following byproduct of the above proof:

\begin{cor}\label{cor:localopen}
With notation as in the lemma, fix a place $w$ of $L$ dividing $p$. The image of the natural restriction map $\Gal(L_w(\mu_{p^\infty})|L_w) \to \Gal(L(\mu_{p^\infty})|L)$ is an open subgroup.
\end{cor}

\begin{proof}
In the course of the above proof we have seen that the image of the restriction map $\Gal(L_w(\mu_{p^\infty})|L_w) \to \Gal(\Q(\mu_{p^\infty})|\Q)$ is an open subgroup of index $\leq B$, but this map factors through $\Gal(L(\mu_{p^\infty})|L)$.   
\end{proof}

\begin{cor}\label{cor:cyclo} With notation as in the lemma, the largest subextension $L_p|L$ of $L^{\rm cyc} | L$ unramified outside the primes dividing $p$ and $\infty$ is $L'(\mu_{p^\infty})$, the extension obtained by adjoining all $p$-power order roots of unity to $L'$. In particular, $L_p$ is a $p$-cyclotomic extension of a field having bounded local degrees.
\end{cor}

\begin{proof}
    The proof is identical to that of \cite[Lemma on p.~316]{Ribet1981CyclotomicTorsion}, except that Ribet works with a finite extension $k|\Q$ and takes the maximal subextension $k'|k$ of $k^{\rm cyc}|k$ unramified at the finite places, whereas we start with our $L$ which may be an infinite extension of $\Q$ but has bounded local degrees.
\end{proof}

\section{Proof of Theorem \ref{mainthm}}\label{secmainthm}

We now return to the situation of the introduction. In particular, $X$ is a smooth proper geometrically integral variety over a number field $K$ with base change $\overline{X}$ to the algebraic closure $\Kbar$. 
We take an  algebraic extension $L|K$ with  maximal cyclotomic extension $L^{\rm cyc}$, all contained in $\overline K$. As before, we write $G_L$ for $\Gal(\Kbar | L)$ and $G_{L^{\rm cyc}}$ for $\Gal(\Kbar | L^{\rm cyc})$.

The remainder of the proof of Theorem \ref{mainthm} is a  generalization of arguments from \cite{RoesslerSzamuely}. We state the two key propositions separately. 

\begin{prop}\label{prop:31-general} Assume that the extension $L|K$ has bounded local degrees.
If $i$ is odd, then
$\left( H^i_{\textup{\'et}}(\overline{X},\mathbf{Q}_p(j))^{\operatorname{ss}} \right)^{G_{L^{\rm cyc}}}=0$ for all primes $p$.
\end{prop}
\begin{proof}
The proof of \cite[Proposition~3.1]{RoesslerSzamuely} goes through almost verbatim with Corollary~\ref{cor:cyclo} replacing Ribet's Galois-theoretic lemma cited above. We briefly review the argument. Note first that we are free to replace $L$ during the argument by algebraic extensions that still have bounded local degrees; in particular, by finite extensions. Pick a simple nonzero
$G_L$-submodule $W\subset \left( H^i_{\textup{\'et}}(\overline{X},\mathbf{Q}_p(j))^{\operatorname{ss}} \right)^{G_{L^{\rm cyc}}}$. Note that $G_L$ acts on $W$ via $G_L/G_{L^{\rm cyc}}$ which is abelian, and that $W$ is semisimple for this action. Using Grothendieck's local monodromy theorem
one then shows that after taking a suitable finite
extension of $L$ the action of $G_L$ on $W$ is unramified outside $p$. By Lemma \ref{lem:bounded local degrees} we are allowed to replace $L$ by the unramified extension $L'|L$ specified by the lemma, so using
Corollary~\ref{cor:cyclo} we may assume that the action of $G_L$ on $W$ factors through $G_p:=\mathrm{Gal}(L(\mu_{p^\infty})|L)$.
For a place $w\mid p$ of $L$ the extension $L_w|\mathbf{Q}_p$ is finite by our assumption on $L$. Moreover, by Corollary \ref{cor:localopen} the image of the restriction map
$\Gal(L_w(\mu_{p^\infty})|L_w)\to G_p$ is open, fixing a finite extension $L''|L$. Replacing $L$ by $L''$
we may thus assume that the local Galois group $\Gal(L_w(\mu_{p^\infty})|L_w)$ surjects onto $G_p$.
As in the original proof, after taking a further finite field extension the existence of the Hodge–Tate decomposition and the Weil conjectures (Deligne's theorem) force an (arithmetic) Frobenius at a place of good reduction $w\nmid p$ to have eigenvalues that are both integral powers of $Nw$ and of absolute value $(Nw)^{j-i/2}$, which is impossible for odd $i$. 
\end{proof}

\begin{prop}\label{prop:34-general} 
Assume that $L|K$ is an abelian extension and $\Gal(L^{\operatorname{cyc}}|\Kcyc)$ has finite exponent. If $i$ is odd, then $\left(H^i_{\textup{\'et}}(\overline{X},\mathbf{Z}/p\mathbf{Z}(j))^{\operatorname{ss}}\right)^{G_{L^{\operatorname{cyc}}}}=0$ for all but finitely many primes $p$.
\end{prop}

The following proof is a direct generalization of the proof of \cite[Proposition~3.4]{RoesslerSzamuely} which is the case $L=K$. However, as the referee pointed out, there are some additional complications in the case $L\supsetneq K$. We isolate the main new technical ingredient in the lemma below. 

Let $K_v|\Qp$ be a finite unramified extension, and denote by $I_v\subset G_{K_v}$ the inertia subgroup. Let moreover $W_v$ be an irreducible $\F_p$-representation of $I_v$ of finite dimension $d$ over $\F_p$. The facts recalled after Definition \ref{def: restricted representations} apply to $W_v$: it is described by a character $\lambda_v:\,I_v\to\F_{p^d}^\times$ that factors through the tame quotient $I_v^t$. The character $\lambda_v$ decomposes as a product of powers of Serre's fundamental characters  $\theta_0,\dots, \theta_{d-1}$ of level $d$; their exponents are the tame inertia weights.

\begin{lem}\label{lem: 34-general}
In the above situation suppose that the tame inertia weights on $W_v$ are contained in an interval $[0,j]$ for some $j>0$. 

If for some integer $0<c<\frac{p-1}{j}$ the character $\lambda_v^c$ equals a power $\bar\chi_p^{r_p}$ of the mod $p$ cyclotomic character $\bar\chi_p$ with  $0 \leq r_p < p-1$, then $c$ divides $r_p$. Moreover, setting $m_p := r_p/c$ we have $0 \leq m_p \leq j$ and $\lambda_v = \bar{\chi}_p^{m_p}$.
\end{lem}

\begin{proof}  We may write $\lambda_v = \prod_{l=0}^{d-1}\theta_{l}^{a_{l}}$, where by assumption the weights $a_{l}$ are integers satisfying $0\leq a_{l}\leq j$ for every $l$.
    
    By \cite[Proposition 8 in \S 1.8]{MR387283}, the mod $p$ cyclotomic character $\overline{\chi}_p$ coincides with the unique fundamental character of level $1$ (recall that $K_v|\mathbf{Q}_p$ is supposed to be unramified). The latter character is the norm of any fundamental character of level $d$, that is, $\overline{\chi}_p = \prod_{l=0}^{d-1}\theta_{l}$. The equality $\lambda_v^c = \overline{\chi}_p^{r_p}$ then gives
    \[
    \prod_{l=0}^{d-1}\theta_{l}^{c a_{l}} = \prod_{l=0}^{d-1}\theta_{l}^{r_p}.
    \]
    Since $0\leq c a_{l}\leq cj<p-1$ and $0 \leq r_p < p-1$, uniqueness of tame inertia weights gives $c a_{l}=r_p$ for all $l$. In particular, $r_p$ is divisible by $c$, and we can write $r_p=c m_p$, where $m_p$ is the common value of all the weights $a_{l}$. Thus $\lambda_v = \prod_{l=0}^{d-1}\theta_{l}^{m_p}=\bar\chi_p^{m_p}$.
\end{proof}

\begin{proof}[Proof of Proposition \ref{prop:34-general}] Throughout the proof, we are free to replace the base field $K$ by a finite extension. Moreover, since $L^{\operatorname{cyc}}$ contains all roots of unity, the statement does not depend on the twist $j$, so we may and do assume that $j\geq i$. As in the proof of \cite[Proposition~3.4]{RoesslerSzamuely}, we may assume that for $p$
large enough $p$ is unramified in $K$, the extension $K(\mu_p)|K$ is nontrivial, the variety $X$ has good reduction at the primes dividing $p$ and we have an isomorphism  $H^i_{\textup{\'et}}(\overline{X},\mathbf{Z}_p(j))/p \, H^i_{\textup{\'et}}(\overline{X},\mathbf{Z}_p(j)) \simeq H^i_{\textup{\'et}}(\overline{X},\mathbf{Z}/p\mathbf{Z}(j))$. Moreover (with an eye towards the lemma above), we assume that the inequality $p> cj+1$ holds, where $c$ is the exponent of the abelian group $\Gal(L^{\operatorname{cyc}}|\Kcyc)$.

Assume now for contradiction that there are infinitely many primes $p$ for which the group $\left(H^i_{\textup{\'et}}(\overline{X},\mathbf{Z}/p\mathbf{Z}(j))^{\operatorname{ss}}\right)^{G_{L^{\operatorname{cyc}}}}$ is nontrivial, and for each such $p$ choose a nonzero irreducible $G_K$-submodule $W_p\subset\left(H^i_{\textup{\'et}}(\overline{X},\mathbf{Z}/p\mathbf{Z}(j))^{\operatorname{ss}}\right)^{G_{L^{\operatorname{cyc}}}}$. The action of $G_K$ on $W_p$ factors via a homomorphism
    $
    \lambda_p : \Gal(L^{\operatorname{cyc}}|K) \to \operatorname{Aut}(W_p).
    $
    Since $\Gal(L^{\operatorname{cyc}}|\Kcyc)$ is abelian of exponent $c$, the homomorphism $\lambda_p^c$ is trivial on $G_{\Kcyc}$ (note that the group $G_K/G_{L^{\operatorname{cyc}}}$ is abelian, so $\lambda_p^c$ is again a homomorphism), so it factors via $\Gal(\Kcyc | K)$. 
Arguing as in \cite[Proposition 3.4]{RoesslerSzamuely}, we find that, after replacing $K$ by a finite extension that can be chosen at the beginning of the proof, the character $\lambda_p^c$ factors through the nontrivial quotient $\Gal(K(\mu_p)|K)$ of $\Gal(\Kcyc | K)$ and hence equals a power $\overline{\chi}_p^{r_p}$ of the mod $p$ cyclotomic character $\overline{\chi}_p$, with  $0 \leq r_p<p-1$.

Since the action of $G_K$ on $W_p$ factors through an abelian group, the image of $G_K$ in
$\Aut_{\F_p}(W_p)$ is contained in
$\Aut_{G_K}(W_p)$. On the other hand, since $W_p$ is irreducible, by Schur's lemma
$\End_{G_K}(W_p)$ is a finite field $\F$ of characteristic $p$.
Thus $\lambda_p$ may be viewed as a character taking values in 
$\F^\times$. 
As such, its image has order prime to $p$.
By Maschke's theorem it then follows that the $G_K$-module $W_p$ remains semisimple upon restriction to any subgroup of $G_K$, in particular to the inertia subgroup $I_v$ of a place $v$ of $K$ above $p$. Pick an irreducible factor $W_v$ of the $I_v$-module $W_p$, and denote by $d$ its dimension over $\F_p$.  As explained at the end of the proof of Theorem \ref{thm: reduction}, in view of our assumption $j\geq i$ Caruso's theorem \cite{Caruso} applies and shows that the tame inertia weights on $W_v$ are contained in the interval $[0,j]$. Now the character $\lambda_v : I_v \to \F_{p^d}^\times$ describing the action of $I_v$ on $W_v$ satisfies $\lambda_v^c = \overline{\chi}_p^{r_p}|_{I_v}$ by the previous paragraph, so since we assumed $p> cj+1$, Lemma \ref{lem: 34-general} applies and gives an equality $r_p=cm_p$ for some integer $m_p\in [0, j]$.

Returning to the equality $\lambda_p^c=\overline{\chi}_p^{r_p}$ of characters  $\Gal(L^{\operatorname{cyc}} | K)\to \F^\times$, we see that the character $\varepsilon_p:=\lambda_p \overline{\chi}_p^{-m_p}: \Gal(L^{\operatorname{cyc}}|K) \to \F^\times$ satisfies $$\varepsilon_p^c= \left( \lambda_p \overline{\chi}_p^{-m_p} \right)^c = \lambda_p^c \overline{\chi}_p^{-r_p} =1,$$ so it takes values in the $c$-th roots of unity in $\F^\times$. 
Since $m_p$ takes values in the finite set $\{0,\ldots,j\}$, after passing to an infinite subset of primes $p$ we may assume that $m_p=m$ is constant. 

Choose now a finite place $w$ of $K$ where $X$ has good reduction. Discarding one more prime if necessary, we may assume that all primes $p$ under consideration are different from the residue characteristic of $w$. Thus both $\lambda_p$ and $\overline{\chi}_p$ are unramified at $w$, and hence so is $\varepsilon_p$.
 Moreover, we have an equality 
    \[
    \lambda_p(\operatorname{Frob}_w) = \varepsilon_p(\operatorname{Frob}_w)\overline{\chi}_p(\operatorname{Frob}_w)^{m} = \varepsilon_p(\operatorname{Frob}_w)(Nw \bmod p)^{m}
    \]
in $\F$, where we know that $\varepsilon_p(\operatorname{Frob}_w)$ is a $c$-th root of unity in $\F^\times$. Since $p>c$ by assumption, we may lift $\varepsilon_p(\operatorname{Frob}_w)$ to a $c$-th root of unity $\zeta\in\overline \Q$. Since there are only finitely many $c$-th roots of unity in $\overline\Q$, we can thus select infinitely many $p$ among our primes for which $\lambda_p(\operatorname{Frob}_w)$ equals the image of $\zeta(Nw)^{m}$ in $\F^\times$ with the {\em same} $\zeta\in \overline\Q$. Let now $P_w(T)$ be the characteristic polynomial of Frobenius at $w$ acting on $H^i_{\textup{\'et}}(\overline{X},\mathbf{Z}_p(j))$. Recall that this polynomial has integral coefficients and is independent of the $p$ under consideration. Since $\lambda_p$ occurs in (the semisimplification of) the mod $p$ reduction of  $H^i_{\textup{\'et}}(\overline{X},\mathbf{Z}_p(j))$, we have $P_w\left(\lambda_p(\operatorname{Frob}_w)\right)=0$  in $\F$.
We conclude that the algebraic integer $P_w(\zeta (Nw)^m)$ has positive valuation for infinitely many places dividing different rational primes $p$, and therefore $P_w(\zeta (Nw)^m)=0$. For odd $i$ this gives a contradiction with Deligne's theorem like at the end of the proof of Proposition \ref{prop:31-general}.
\end{proof}

\begin{remark}
As noted above, the above proof reduces to that of \cite[Proposition~3.4]{RoesslerSzamuely} in the case $L=K$. We take this opportunity to point out that the important reduction to the case $j\geq i$ was unfortunately missing from {\em loc. cit.}     
\end{remark}

Finally, we arrive at:

\begin{proof}[Proof of Theorem \ref{mainthm}]
Suppose by contradiction that $H^i_{\textup{ét}}(\overline{X}, \Q/\Z(j))^{G_{K^{\rm Kum}}}$ is infinite. By Theorem \ref{thm: reduction}, there exist a finite extension $K' | K$ and an abelian subextension $F$ of $K' \Kkum | K'$ such that at least one of the following holds:
\begin{enumerate}
    \item $\left( H^i_{\textup{ét}}(\overline{X}, \Q_p(j))^{\operatorname{ss}} \right)^{G_F} \neq (0)$ for some prime $p$;
    \item $\left( H^i_{\textup{ét}}(\overline{X}, \Z/p\Z(j))^{\operatorname{ss}} \right)^{G_F} \neq (0)$ for infinitely many primes $p$.
\end{enumerate}
Replacing $K$ with $K'$ (which is again a number field), we have an abelian extension $F|K$ contained in $K^{\operatorname{Kum}}$ for which (1) or (2) above holds. Applying Proposition \ref{prop: abelianisation Kummer}, we find that $F$ is  contained in the compositum of $K^{\operatorname{cyc}}$ and an abelian extension $L|K$ of finite exponent. In particular,  $G_{L^{\operatorname{cyc}}}$ is contained in $G_F$ and therefore (1) or (2) above holds with  $G_F$ replaced by $G_{L^{\operatorname{cyc}}}$, contradicting Proposition \ref{prop:31-general} or Proposition \ref{prop:34-general}.
\end{proof}

\begin{remark}
    
 The finiteness of $H^i_{\textup{ét}}(\overline{X}, \Q_p(j))^{G_{K^{\rm Kum}}}$ for odd $i$ was also obtained in \cite[Theorem 2.11]{MurotaniOzeki} using results from \cite{RoesslerSzamuely}.
\end{remark}

\section{Proof of Proposition \ref{prop:solvable_example}}\label{secexample}

In this section we describe the example announced in Proposition \ref{prop:solvable_example}. We begin with the construction of particular infinite solvable extensions of $\Q$. 

\begin{prop}\label{prop: real quadratic construction}
Let $K$ be a real quadratic field of odd conductor and  class number $1$. 
There exists an infinite sequence of rational primes $ p_1 < p_2 < \ldots$, each split in $\mathbb{Q}(i)$, and a totally real Galois extension $M|\Q$ with $\Gal(M | K)\cong H $, where

$$H := \prod_{i \geq 1} H_{i}$$
with
\[
H_{i} := (\Z/p_i\Z)^\times \times (\Z/p_i\Z)^\times,
\]
such that moreover the group extension
$$1\to H\to \Gal(M | \Q)\to \Z/2\Z\to 1$$
corresponding to the tower of fields $M|K|\Q$
splits as a semidirect product
$$\Gal(M | \Q)\cong H \rtimes \Z/2\Z,$$ with the nontrivial element of $\Z/2\Z$ acting on $H$  by swapping the $(\Z/p_i\Z)^\times$-factors of each $H_i$.  
\end{prop}

\begin{proof}
We will obtain $M$ as the increasing union of certain abelian extensions $L^{(n)}$ that satisfy
\[
\Gal(L^{(n)} | K) \cong \prod_{i=1}^n H_{i},
\]
with $\Z/2\Z$ acting on $\prod_{i=1}^n H_{i}$ as in the statement. In turn, the extension $L^{(n)}$ will be constructed as the compositum of extensions $L_{p_1}, \ldots, L_{p_n}$, so we start by describing the construction of a single $L_{p}$.

Let $c$ be the conductor of $K$. We will choose our primes to be congruent to $1$ modulo $c$ and to $5 \pmod 8$. These conditions ensure that
\begin{enumerate}
    \item $p$ splits in $\mathcal{O}_K$ (the condition $p \equiv 1 \pmod{c}$ implies that $p$ splits in $\Q(\zeta_c) \supseteq K$);
    \item $4 \mid p-1$ but $8 \nmid p-1$. In particular, we will have $p \equiv 1 \pmod{4}$, so that $p$ splits in $\mathbb{Q}(i)$.
\end{enumerate}
Since $c$ is odd by assumption, the conditions are compatible and hence Dirichlet's theorem implies that there exist infinitely many such primes. Fix one such $p$ and write $p\mathcal O_K = \mathfrak{p}\mathfrak{p}'$. By the Grunwald--Wang theorem (see \cite[Chapter X, Theorem 5]{ArtinTate} or \cite[Theorem 9.2.8]{nsw}) there is a cyclic extension $M_p|K$ with group $(\Z/p\Z)^\times$ such that:
\begin{itemize}
\item there is only one prime in $\mathcal O_{M_p}$ above $\mathfrak{p}$ and the corresponding local extension is unramified and cyclic with group $(\Z/p\Z)^\times$ (so the
      decomposition group at $\mathfrak{p}$ is all of $\Gal(M_p | K)$);
\item $\mathfrak{p}'$  splits completely in $M_p$ (so the corresponding local extensions are all trivial);
\item $M_p$ is a totally real number field (i.e., the infinite places of $K$ split in $M_p$).
\end{itemize}

Note that condition (2) implies that we are not in the so-called \textit{special case} of the Grunwald--Wang theorem, because in the notation of \cite[Chapter X, Theorems 1 and 5]{ArtinTate} we have $m=2^2 m'$ with $m'$ odd and the inequality $2>s \geq 2$ is not satisfied. So the theorem indeed applies and gives an extension $M_p|K$ with the local behaviour described above.

The extension $M_p|K$ is defined by a character $\chi$ of $\Gal(K^{\rm ab} |  K)$. The conjugation action of the nontrivial element $\sigma\in\Gal(K|\Q)$ on $\Gal(K^{\rm ab}|K)$ coming from the extension
$$
1\to \Gal(K^{\rm ab}|K)\to \Gal(K^{\rm ab}|\Q)\to \Gal(K|\Q)\to 1
$$
defines another character $\chi^\sigma$ whose kernel fixes a Galois extension $M_p^\sigma|K$ with group $(\Z/p\Z)^\times$. Its local behaviour at the primes above $p$ is obtained by swapping
$\mathfrak{p}$ and $\mathfrak{p}'$: the prime $\mathfrak{p}'$ is totally inert in $M_p^{\sigma}|K$ and $\mathfrak{p}$ is completely split.  Therefore the cyclic extensions $M_p$ and $M_p^\sigma$ are linearly disjoint over $K$ as $\mathfrak{p}$ is totally inert in the first and splits completely in the second. 
Let $L_p := M_p M_p^\sigma$ be their compositum (it is also the Galois closure of $M_p$ over $\Q$). We have $\Gal(L_p|K)\cong (\Z/p\Z)^\times\times (\Z/p\Z)^\times$; moreover, $\sigma$ acts on $\Gal(L_p|K)$ by swapping the factors. With this action of $\Z/2\Z$ the abelian group $(\Z/p\Z)^\times\times (\Z/p\Z)^\times$ is none but the induced $\Z/2\Z$-module ${\rm Ind}^{\Z/2\Z}_{\{1\}}((\Z/p\Z)^\times)$ and therefore the cohomology group $H^2(\Z/2\Z,(\Z/p\Z)^\times\times (\Z/p\Z)^\times)$ vanishes. This implies that $\Gal(L_p|\Q)$ is a split extension of $\Gal(K|\Q)$ by $\Gal(L_p|K)$, i.e.~a semidirect product $((\Z/p\Z)^\times\times (\Z/p\Z)^\times)\rtimes \Z/2\Z$.

We now proceed inductively. Assume $L_{p_1},\dots, L_{p_n}$ have been constructed, and let $p_{n+1}$ be a prime strictly larger than $p_1, \ldots, p_{n}$, congruent to $1$ modulo $c$ and to $5$ modulo $8$, and unramified in each of the fields $L_{p_1}, \ldots, L_{p_n}$. As noted above, there are infinitely many primes that satisfy the congruence conditions, and the unramifiedness condition  excludes only finitely many of them. So such a $p_{n+1}$ indeed exists. 
We carry out the above construction to obtain a field $L_{p_{n+1}}|K$ with
\[
\Gal(L_{p_{n+1}} | K)\cong H_{{n+1}}=(\Z/p_{n+1}\Z)^\times\times (\Z/p_{n+1}\Z)^\times.
\]
When we apply the Grunwald--Wang theorem as explained above, we further impose that the extension $L_{p_{n+1}}|K$ be locally trivial (hence in particular unramified) at all the finitely many primes above any of the rational primes that ramify in $L_{p_1}, \ldots, L_{p_n}$.
It follows that the fields $L_{p_i}$ for $1 \leq i \leq n+1$ are linearly disjoint over $K$. Indeed, the above ramification condition implies that the intersection of one of them with the compositum of the others is an everywhere unramified abelian extension of $K$. Since $K$ has class number $1$, this implies that the intersection must be $K$ itself. Setting $L^{(n+1)} := L_{p_1}\cdots L_{p_{n+1}}$, we then have 
\[
\Gal(L^{(n+1)} | K)\;\cong\;\prod_{i=1}^{n+1} \Gal(L_{p_i} | K)
             \;\cong\;\prod_{i=1}^{n+1} H_{i}.
\]
By construction we have $\Gal(L_{p_i}| \Q)\cong H_i\rtimes \Z/2\Z$ with the desired action of $\sigma$ on $ H_{i}$ for each $i$, and moreover $L^{(n+1)}$ is stable by the action of  $\Z/2\Z$. Hence $L^{(n+1)}$ is a Galois extension of $\Q$ which is moreover totally real as so are the $L_{p_i}$. 

Finally, as already noted at the beginning of the proof, to complete the construction it suffices to take $M=\bigcup_n L^{(n)}$. 
\end{proof}

\begin{remark}
Related results on the existence of Galois extensions of $\Q$ with solvable Galois group and prescribed local conditions are proven, for example, in \cite[Theorem 6.6]{Neukirch}. However, it seems that Neukirch's statements do not immediately apply in the context of our inductive construction.
\end{remark}

We now construct the other field extension required by Proposition \ref{prop:solvable_example}.

\begin{prop}\label{prop:solvable_example2} Let $E$ be an elliptic curve  defined over $\Q$ that has CM by $\Z[i]$, and let $p_1<p_2< \ldots$ be a sequence of primes of good reduction for $E$ that split completely in $\Q(i)$. There exists a Galois extension $M'|\Q$ containing $\Q(i)$ such that $E$  has infinite torsion over $M'$ and moreover 
$\Gal(M' | \Q)\cong H\rtimes \Z/2\Z$, where $H= \prod_i((\Z/p_i\Z)^\times \times (\Z/p_i\Z)^\times )$ and $\Z/2\Z$ acts by swapping the $(\Z/p_i\Z)^\times$-factors.     
\end{prop}

\begin{proof} 
Let first $p$ be any prime of good reduction for $E$ that splits in $\Q(i)$, and consider the representation $\rho_p:{\operatorname{Gal}(\overline{\Q} | \Q(i))}\to \Aut_{\F_p}(E[p])$ on the $p$-torsion points of an elliptic curve $E$ as above. Since $E$ has CM by $\Z[i]$, Corollary 2 to Theorem 5 of \cite{serretate} implies that $\rho_p$ factors through an injection  $(\F_p^2)^\times\hookrightarrow \End_{\F_p}(E[p])$, yielding a map ${\operatorname{Gal}(\overline{\Q} | \Q(i))}^{\operatorname{ab}}\to (\F_p^2)^\times$.  On the other hand, since $\Q(i)$ has class number $1$, class field theory gives a  surjection $r:\,\prod_{v \in \Omega_{\Q(i)}} \mathcal{O}_{\Q(i), v}^\times \to {\Gal(\overline{\Q} | \Q(i))^{\operatorname{ab}}}$ where $\Omega_{\Q(i)}$ denotes the set of finite places of $\Q(i)$. Composing the two we obtain a map $\prod_{v \in \Omega_{\Q(i)}} \mathcal{O}_{\Q(i), v}^\times\to (\F_p^2)^\times$. By the results of \cite{serretate} summarized in  \cite[Theorem 5.1]{MR3766118}, this map sends $(a_v)_{v \in \Omega_{\Q(i)}}$ to $(a_{\mathfrak{p}}^{-1}, a_{\mathfrak{p'}}^{-1})$ mod $p$. 
(Indeed, note that in the formula of \textit{loc.~cit.} the maps $N_{K_\ell/E_\ell^*}$,  $\Phi_{E,S}$ and $\varepsilon$ are all trivial in our case, the first two because we are dealing with a CM elliptic curve over $\Q(i)$ and the last one because $\mathfrak{p}$, $\mathfrak{p}'$ are primes of good reduction for $E$.) In particular, this map is the sum of two surjective characters with values in $(\Z/p\Z)^\times$.

Now let $p_1 < p_2 < \ldots$ be the sequence of completely split primes of good reduction as in the assumption. With this choice of primes consider the product representation $$\prod_i\rho_{p_i}:\, {\operatorname{Gal}(\overline{\Q} | \Q(i))}\to \prod_i\Aut_{\F_{p_i}}(E[p_i]).$$ Factoring through the $(\F_{p_i}^2)^\times$ and composing with $r$ as above we obtain a surjection
\[
 \prod_{i} \left( \mathcal{O}_{\Q(i), \mathfrak{p}_{i}}^\times \times \mathcal{O}_{\Q(i), \mathfrak{p}'_{i}}^\times \right) \twoheadrightarrow \prod_{i} ((\Z/p_i\Z)^\times \times (\Z/p_i\Z)^\times)=H.
 \]
 Thus $H$ arises as a quotient of  $\operatorname{Gal}(\overline{\Q} | \Q(i))$, corresponding to a field extension $M'|\Q(i)$. Since the action of $\Gal(\Q(i) | \Q)$ exchanges the primes $\mathfrak{p}_i$ and $\mathfrak{p}_i'$, its action on $\Gal(M'| \Q(i))$ exchanges the $(\Z/p_i\Z)^\times$-factors in $H$. {Thus $H$ is again an induced $\Z/2\Z$-module, so as in the previous proof we obtain that} $\Gal(M'|\Q)$ is a semidirect product $H\rtimes \Z/2\Z$ as in the statement.

Finally, note that by construction $\operatorname{Gal}(\overline{\Q}| M')$ fixes the $p_i$-torsion points of $E$ for all $i$, so that $E$ has infinite torsion over $M'$.
\end{proof}

\begin{proof}[Proof of Proposition \ref{prop:solvable_example}] Take an elliptic curve $E$ as in Proposition \ref{prop:solvable_example2} and choose a sequence of primes $p_1<p_2<\cdots$ as in Proposition \ref{prop: real quadratic construction}. Leaving out finitely many of the $p_i$ we may assume they are all of good reduction for $E$. Construct the extensions $M$ and $M'$ of Propositions \ref{prop: real quadratic construction} and \ref{prop:solvable_example2} with this sequence of primes. They have isomorphic Galois groups over $\Q$ that are solvable of class 2, and we also know that $E$ has infinite torsion over $M'$. On the other hand,  note that every abelian variety defined over the totally real field $M$ is defined over some totally real subfield $\widetilde M\subset M$ finite over $\Q$ and in particular cannot have CM defined over $\widetilde M$. {Without loss of generality, we can assume that $\widetilde{M}$ contains $K$, so that $M|\widetilde{M}$ is an abelian extension (being a subextension of $M|K$).}
Therefore by the result of Zarhin already cited in the introduction (\cite{ZarhinCM}, Corollary of Theorem 1) it has finite torsion over every abelian extension of $\widetilde{M}$, and in particular over $M$.
\end{proof}

\begin{remark}
To be somewhat more concrete, in Proposition \ref{prop:solvable_example2} we can work with the elliptic curve $E$ of Weierstrass equation $y^2=x^3+x$ that has CM by $\Z[i]$ and 2 as its only prime of bad reduction. Then in Proposition \ref{prop: real quadratic construction} we may take $K=\Q(\sqrt{5})$ that has class number 1 and conductor 5. For the sequence of primes $p_1< p_2<\ldots$ we can take any sequence provided by Proposition \ref{prop: real quadratic construction} provided that $p_1\neq 2$.     
\end{remark}

\bigskip

\noindent{\bf Funding}.
This work was supported by the Italian \textit{Ministero dell'Università e della Ricerca} through grants PRIN-2022HPSNCR and PRIN-2022BTA242 (funded by the European Union project Next Generation EU).

\bibliographystyle{plain}
\bibliography{biblio}
\end{document}